\documentclass[preprint]{elsarticle}
           \usepackage[T2A]{fontenc}   
           \usepackage[utf8]{inputenc}
           \usepackage[english]{babel}
\usepackage{amsmath}
\usepackage{amsfonts}
\usepackage{amssymb}
\usepackage{amsthm}
\usepackage[pdftex,unicode]{hyperref}

\newtheorem{assertion}{Statement}

\newtheorem{proposition}{Proposition}
\newtheorem{theorem}{Theorem}
\newtheorem*{theorem*}{Theorem}

\newtheorem*{lemma*}{Lemma}
\newtheorem*{claim}{Claim}
\newtheorem{cor}{Corollary}

\theoremstyle{definition}
\newtheorem{example}{Example}
\newtheorem{definition}{Definition}
\newtheorem{question}{Question}

\theoremstyle{remark}

\newtheorem*{note*}{Remark}

\def\sset#1{\{#1\}}

\def\set#1{\bbset#1\eeset}
\def\bbset#1:#2\eeset{\{#1\,:\,#2\}}

\def\bbsett#1:#2\eesett{\{#1\,:\,\text{#2}\}}

\def\si{\sigma}
\def\om{\omega}
\def\term#1{{\it #1}}
\def\nne{\cN_e}
\def\cN{\mathcal{N}}
\def\gi{{\mathfrak{i}}}
\def\gm{{\mathfrak{m}}}
\def\ga{\gamma}
\def\la{\lambda}
\def\ph{\varphi}
\def\de{\delta}
\def\cl#1{\overline{#1}}
\def\Int{\operatorname{Int}}
\def\cS{{\mathcal S}}
\def\tp{\Tau}

\def\lrarr{\Leftrightarrow}
\def\al{\alpha}
\def\es{\varnothing}
\def\Tau{{\mathcal T}}
\def\cP{{\mathcal P}}
\def\cF{{\mathcal F}}
\def\be{\beta}
\def\sq#1#2{(#1)_{#2}}
\def\sqn#1{\sq{#1}{n\in\om}}
\def\sqnn#1{\sqn{#1_n}}
\def\nom{{n\in\om}}
\def\id{\operatorname{id}}

\def\clp#1{\cl{\sset{#1}}}
\def\tsf#1{{#1}_{T_0}}
\def\pitsf{{\pi}_{T_0}}
\def\PT{$\mathrm{(\star)}$}
\def\simt{\mathbin{{\sim_T}}}
\def\gj{{\mathfrak{j}}}

\def\sym{\operatorname{Sym}}
\def\pUs/{$\mathrm{(U_s)}$}
\def\Si{\Sigma}
\def\cC{\mathcal{C}}

\begin{document}

\begin{frontmatter}

\title{Almost paratopological groups}
\author{Evgenii Reznichenko}

\ead{erezn@inbox.ru}

\address{Department of General Topology and Geometry, Mechanics and Mathematics Faculty, M.~V.~Lomonosov Moscow State University, Leninskie Gory 1, Moscow, 199991 Russia}

\begin{abstract}
A class of almost paratopological groups is introduced, which (1) contains paratopological groups and Hausdorff quasitopological groups; (2) is closed under products; (3) subgroups.
Almost paratopological $T_1$ groups $G$ are characterized by the fact that $\set{(x,y)\in G^2: xy=e}$ is closed in $G^2$.
A compact almost paratopological group is topological.
A regular $\Si$-space with countable extend and a separately continuous Mal’tsev operation is $\om$-cellular (and  ccc).
A $\si$-compact regular almost paratopological group is ccc.
In particular, a $\si$-compact regular quasitopological group is ccc.
\end{abstract}
\end{frontmatter}
\section{Introduction}\label{sec:intro}
One of the most important concepts in mathematics is the concept of a topological group.
A topological group is a group with a topology with respect to which the operations of product and inverse are continuous.
Other groups with topology are also important. Groups $G$ with topology is called
\begin{itemize}
\item
\term{paratopological} if multiplication in $G$ is continuous;
\item
\term{semitopological} if multiplication in $G$ is separately continuous;
\item
\term{quasitopological} if $G$ is a semitopological group and the operation
taking the inverse element $x\mapsto x^{-1}$ is continuous.
\end{itemize}
In the book \cite{at2009} and the survey \cite{tk2014} one can find further information about the enumerated classes of groups.
This article introduces the concept of almost paratopological group, a class of semitopological groups that contains paratopological groups and Hausdorff quasitopological groups. This class of semitopological groups is closed under the product and subgroups (Theorem \ref{t:apt:1}). For them, assertions are proved that were previously known for paratopological groups.
A compact almost paratopological group is topological (Theorem \ref{t:apt:2}).
In particular, a compact paratopological group is topological \cite[Lemma 6]{Ravsky2019}.

A mapping $M:X^3\to X$ is called a \term{Mal'tsev operation} if 
\[
M(x,y,y)=M(y,y,x)=x
\]
for all $x,y\in X$. There is a natural Mal'tsev operation on groups:
\[
M_G(x,y,z)=xy^{-1}z
\]
for $x,y,z\in G$. On retracts of sets with the Mal'tsev operation, there is a Mal'tsev operation: if $r: X\to X$ is a retraction, then we set $M'(x,y,z)=r(M(x,y,z))$ for $x, y,z\in r(X)$.
A topological space with a continuous Mal'tsev operation is called a \term{Mal'tsev spaces}.
On topological groups, the Mal'tsev operation $M_G$ is continuous, and on quasitopological groups $M_G$ is separately continuous.

In \cite{ReznichenkoUspenskij1998} it was proved that if $X$ is a Lindelöf $\Si$ or a countably compact Tychonoff space with a separately continuous Mal'tsev operation, then $X$ is an $\om$-cellular space (in particular, $X$ is ccc). It was also noted there that this assertion can be proved for Tychonoff $\Si$-spaces with  countable extend. This article proves that if $X$ is a regular $\Si$-space with  countable extend with a separately continuous Mal'tsev operation, then $X$ is an $\om$-cellular space (Theorem \ref{t:maltcev:cel:1}) .

The above facts imply that a regular quasitopological $\Si$ group with  countable extend is $\om$-cellular.
Based on this fact, it is proved that the Lindelöf $\Si$ almost paratopological group is $\om$-cellular (Theorem \ref{t:ccc:3}).
In particular, an $\si$-compact regular quasitopological or almost paratopological group is ccc.

\section{Definitions and notation}\label{sec:defs}
The identity in the group $G$ will be denoted as $e$, the family of open neighborhoods $e$ as $\nne$.
Denote the multiplication in the group
\[
\gm: G\times G\to G,\ (x,y)\mapsto xy
\]
and operation
taking the inverse element
\[
\gi: G\to G,\ x\mapsto x^{-1}.
\]

Let $X$ be a topological space. Then two points $x$ and $y$ in $X$ are \term{topologically indistinguishable} if they have exactly the same neighborhoods; that is, $x\in\clp y$ and $y\in\clp x$.
Points $x$ and $y$ in $X$ are \term{topologically distinguishable} if they are not topologically indistinguishable.

The space $X$ has \term{countable extend} if every discrete closed subset is at most countable.

A space $X$ is \term{ccc}, or has the \term{Suslin property},
if every disjoint family of nonempty open sets in $X$ is at most countable.
Let us say that a space $X$ is \term{$\om$-cellular}
if for any family $\la$ of $G_\de$-subsets of $X$ there exists a countable subfamily $\mu\subset \la$ such
that $\cl{\bigcup\la}=\cl{\bigcup\mu}$. Clearly $\om$-cellular spaces are ccc.

A space $X$ is called a \term{$\Si$-space} if
there exists a $\si$-locally finite family
$\cS$ and the covering of $\cC$ by closed countably compact sets, such that if $C\in\cC$ and $U\supset C$ is open, then $C\subset S \subset U$ for some $S\in \cS$.
Lindelöf $\Si$-spaces are exactly the regular continuous images of perfect inverse images of metrizable separable spaces.
$\si$-countably compact spaces are $\Si$-spaces.
Regular $\si$-compact spaces are Lindelöf $\Si$-spaces.

\begin{assertion}\label{a:maltcev:cel:1}
A space $X$ is a $\Si$-space with  countable extend if and only if
there is a countable family
$\cS$ and the covering of $\cC$ by closed countably compact sets, such that if $C\in\cC$ and $U\supset C$ are open, then $C\subset S \subset U$ for some $S\in \cS$.
\end{assertion}
\begin{proof}
Let $\cS$ be a $\si$-discrete family and $\cC$ be a covering of closed countably compact sets, such that if $C\in\cC$ and $U\supset C$ are open, then $C\subset S \subset U$ for some $S\in\cS$.
Since $X$ has a countable extend, then $\cS$ is countable.
Conversely, it suffices to note that a countable family is a $\si$-discrete family.
\end{proof}

A space is \term{$\si$-countably compact} if it is the union of a countable family of countably compact subspaces.

Let us call a space \term{closely $\si$-countably compact} if it is the union of a countable family of closed countably compact subspaces.
Any closely $\si$-countably compact space is a $\Si$-space with  countable extend.

\section{$R_0$ spaces and groups}\label{sec:rsa}
Let $X$ and $Y$ be topological spaces. We say that a map $f: X\to Y$ \term{preserves the topology} ($f$ is \term{topology-preserving}) if
\begin{itemize}
\item[\PT] the mapping $f$ is surjective, continuous, open and closed, and the set $U\subset X$ is open if and only if $U=f^{-1}(f(U))$ and $f (U)$ open.
\end{itemize}

It is easy to verify the following assertion.

\begin{proposition}\label{p:rsa:1}
Let $(X,\tp_X)$ and $(Y,\tp_Y)$ be topological spaces, $f: X\to Y$ be a surjective mapping. The following conditions are equivalent:
\begin{enumerate}
\item
$f$ preserves the topology;
\item
mapping
\[
f^{\#}: \tp_Y\to \tp_X,\ U\mapsto f^{-1}(U)
\]
is a bijection;
\item
$f$ is quotient and the subspace $f^{-1}(y)$ is antidiscrete for every $y\in Y$.
\end{enumerate}
\end{proposition}

The relation `$x$ and $y$ are topologically indistinguishable' on $X$ is an equivalence relation, denote this equivalence as $\simt$: $x\simt y$ if and only if $x\in \clp y$ and $ y\in\clp x$.

  The space quotient by this equivalence relation is denoted as $\tsf X$ and $\pitsf: X\to \tsf X$ is the quotient mapping. The space $\tsf X$ is a $T_0$ space and the mapping $f=\pitsf: X\to Y=\tsf X$ preserves the topology.
If $X$ is a $T_0$ space, then $\pitsf$ is a homeomorphism.
Recall the axiom of separation of $R_0$ topological spaces.

The space $X$ is $R_0$, or \term{symmetric}, if
$x \in \clp y$ implies $y \in \clp x$ (that is, $x\simt y$) for any $x,y\in X$ \cite{Schechter1996}.

The space $X$ is $R_0$ if and only if $\tsf X$ is $T_1$. See Section 16.2 in \cite{Schechter1996}.

\subsection{Arkhangel'skii--Motorov theorem}

In 1984 D.B. Motorov constructed a famous example of a compact space that is not a retract of any compact space.
This is the closure on the plane of the graph of the function $\sin \frac 1x$ with domain $(0,1]$ (\cite{Seminar1985}).
A little later, A.V. Arkhangelsky published Motorov's improved results \cite{Arkhangelskii1985}.
We need the following two statements from their papers.

\begin{theorem}[{Arkhangel'skii \cite[Corollary 1]{Arkhangelskii1985}}]\label{t:rsa:mot:1}
Let $X$ be a homogeneous space and  $\clp x$ is compact for $x\in X$. Then $X$  is a $R_0$ space.
\end{theorem}

\begin{theorem}[Motorov]\label{t:rsa:mot:2}
If $X$ is a homogeneous compact space, then $X$ is an $R_0$ space.
\end{theorem}

\subsection{Topology-preserving group homomorphisms}

\begin{proposition}\label{p:rsa:2}
Let $G$ be a semitopological group, $H$ a group with topology, and $\ph: G\to H$ a surjective quotient homomorphism.
Then $\ph$ is open and $H$ is a semitopological group.
\end{proposition}
\begin{proof}
Let $U$ be an open subset of $X$. Because
\[
\ph^{-1}(\ph(U))= U\ker \ph
\]
and right shifts are homeomorphisms, then $\ph^{-1}(\ph(U))$ and hence $\ph(U)$ are open.

Let $V$ be an open subset of $H$, $h\in H$ and $g\in \ph^{-1}(h)$. Because
\begin{align*}
\ph^{-1}(hV)&=g\ph^{-1}(V)
&&\text{and}&
\ph^{-1}(Vh)&=\ph^{-1}(V)g,
\end{align*}
then the sets $hV$ and $Vh$ are open.
\end{proof}

Propositions \ref{p:rsa:1} and \ref{p:rsa:2} imply the following assertion.

\begin{proposition}\label{p:rsa:3}
Let $G$ be a semitopological group, $H$ a group with topology, and $\ph: G\to H$ a surjective homomorphism.
  The following conditions are equivalent:
\begin{enumerate}
\item
$\ph$ preserves the topology;
\item
the mapping $\ph$ is quotient and $\ker \ph$ is antidiscrete.
\end{enumerate}
If $\ph$ preserves the topology, then $H$ is a semitopological group.
\end{proposition}

The following theorem is practically the same as \cite[Theorem 3.1]{Tkachenko2014}.

\begin{theorem}\label{t:rsa:1}
Let $G$ be a semitopological group,
$H = \clp e \cap \clp e ^{-1}$.
Then $H$ is a normal antidiscrete subgroup, the quotient group $G/H$ coincides with $\tsf G$, $G/H$ is a semitopological $T_0$ group, and the quotient mapping $G\to G/H$ is topology-preserving.
\end{theorem}
\begin{proof}
Let $x\in G$. Then $x\in H$ $\lrarr$ $x\in\clp e$ and $x^{-1}\in\clp e$ $\lrarr$ $x\in\clp e$ and $e\in \clp x$ $\lrarr$ $e\simt x$.
Hence $H=\pitsf^{-1}(\pitsf(e))$.

The equivalence relation $\simt$ is invariant under right and left shifts, so the quotient set by this equivalence relation coincides with the right (and left) cosets with respect to the normal subgroup $H$.
It remains to apply Proposition \ref{p:rsa:3}.
\end{proof}

Theorems \ref{t:rsa:1} and \ref{t:rsa:mot:1} imply the following assertion.

\begin{theorem}\label{t:rsa:2}
Let $G$ be a semitopological group, $H=\clp e$ be compact. Then $G$ is the $R_0$ space,
$H$ is a normal closed antidiscrete subgroup, the quotient group $G/H$ coincides with $\tsf G$, $G/H$ is a semitopological $T_1$ group, and the quotient mapping $G\to G/H$ is topology-preserving.
\end{theorem}

\section{Almost paratopological groups}\label{sec:apt}
\begin{definition}
A semitopological group $G$ is called \term{almost paratopological},
if for any $g\in G$ such that $e\notin \clp g$ there exists a neighborhood $U$ of $e$ such that $g\notin U^2$.
\end{definition}

\begin{theorem}\label{t:apt:1}
Let $G$ be a semitopological group. If any of the following conditions is satisfied, then $G$ is almost paratopological.
\begin{enumerate}
\item
$G$ is a subgroup of an almost paratopological group;
\item
$G$ is the product of almost paratopological groups;
\item
$G$ is a paratopological group;
\item
$G$ is a Hausdorff quasitopological group;
\item
there exists a continuous isomorphism of $G$ onto a $T_1$ almost paratopological group.
\end{enumerate}
\end{theorem}
\begin{proof}
Conditions (1) and (3) are obvious.

(2) Let $G=\prod_{\al\in A}G_\al$, where $G_\al$ is an almost paratopological group for all $\al\in A$.
Let $g=(g_\al)_{\al\in A}\in G$ and $e=(e_\al)_{\al\in A}\notin \clp g$. Then $e_\al\notin\clp{g_\al}$ for some $\al\in A$.
There is a neighborhood $U_\al$ of $e_\al$ such that $g_\al\notin U_\al^2$. Let $U$ be the inverse image of $U_\al$ under the projection of $G$ onto $G_\al$. Then $g\notin U^2$.

(4) Let $g\in G\setminus \sset e$. There is $U\in\nne$ for which $U=U^{-1}$ and $U\cap gU=\es$. Then $g\notin U^2$.

(5) Let $\ph: G\to H$ be a continuous isomorphism of the group $G$ onto  a $T_1$ almost paratopological group $H$. Let $g\in G\setminus \sset e$. There is a neighborhood $V\subset H$ of $\ph(e)$ such that $\ph(g)\notin V^2$. Then $g\notin U^2$, where $U=\ph^{-1}(V)$.
\end{proof}

\begin{example}\label{e:apt:1}
Let $G$ be the group of integers with topology consisting of co-finite sets. The group $G$ is a quasitopological compact $T_1$ group which is not almost paratopological.
\end{example}

\begin{proposition}\label{p:apt:1}
Let $G$ be a semitopological group, $M\subset G$. Then
\[
\cl M = \bigcap_{U\in\nne} MU^{-1}=\bigcap_{U\in\nne} U^{-1}M.
\]
\end{proposition}
\begin{proof}
Let $g\in \cl M$ and $U\in\nne$. Then $gU\cap M\neq\es$ and $g\in MU^{-1}$. Hence $\cl M \subset \bigcap_{U\in\nne} MU^{-1}$.
Let $g\in \bigcap_{U\in\nne} MU^{-1}$. Then $gU\cap M$ for any $U\in\nne$. Hence $\bigcap_{U\in\nne} MU^{-1} \subset \cl M$.
$\cl M = \bigcap_{U\in\nne} U^{-1}M$ is proved similarly.
\end{proof}

For a group $G$ we denote
\begin{align*}
S_G &= \set{(x,y)\in G\times G: xy=e},
&
E_G &= \bigcap_{U\in\nne} \cl{U^{-1}}.
\end{align*}

\begin{proposition}\label{p:apt:2}
Let $G$ be a semitopological group.
Then
\begin{enumerate}
\item
\[
\clp e \subset E_G = \bigcap_{U\in\nne} \left(U^{-1}\right)^{2};
\]
\item
the group $G$ is almost paratopological if and only if $E_G=\clp e$;
\item
$\cl{S_G}=\gm^{-1}(E_G)$;
\item
the following conditions are equivalent:
\begin{enumerate}
\item
$G$ is an almost paratopological $T_1$ group;
\item
$E_G=\sset e$;
\item
$S_G$ is closed in $G^2$.
\end{enumerate}
\end{enumerate}
\end{proposition}
\begin{proof}
Denote $Q=\bigcap_{U\in\nne} \left(U^{-1}\right)^{2}$.

(1). Proposition \ref{p:apt:1} implies that $\clp e \subset E_G \subset Q$.
Let $x\in G\setminus E_G$. Then $xU\cap U^{-1}=\es$ for some $U\in\nne$. Hence $x\notin (U^{-1})^2\supset Q$. We get $Q\subset E_G$.

(2). Assume that $G$ is an almost paratopological group. It follows from (1) that $\clp e \subset E_G$.
Let $g\in G\setminus \clp e$. Then $e\notin \clp{g^{-1}}$. Since $G$ is almost paratopological, it follows that $g^{-1}\notin U^2$ for some $U\in\nne$. Then $g\notin (U^{-1})^2\supset E_G$. Hence $E_G\subset \clp e$.

Suppose $E_G=\clp e$. Let $g\in G$ and $e\notin\clp g$. Then $g^{-1}\notin \clp e= E_G$. It follows from (1) that $g^{-1}\notin (U^{-1})^2$ for some $U\in\nne$. We get, $g\notin U^2$.

(3). Let $(x,y)\in G$. Then $(x,y)\in \cl{S_G}$ $\lrarr$ $(Ux \times Uy)\cap S_G$ for all $U\in\nne$ $\lrarr$
$e\in UxUy$ for all $U\in\nne$ $\lrarr$ $e\in UUxy$ for all $U\in\nne$ $\lrarr$ $xy\in (U^{-1})^2$ for all $U\in\nne$ $\lrarr$
$xy\in Q$. It follows from (1) that $(x,y)\in \cl{S_G}$ $\lrarr$ $xy\in E_G$.

(4). From (2) follows (a)$\lrarr$(b). Since $S_G=\gm^{-1}(e)$, then (3) implies (b)$\lrarr$(c).
\end{proof}

For a group with the topology $(G,\tp)$, we denote
\[
\tp_{Sym}=\set{U\cap V^{-1}: U,V\in \tp}
\]
and $\sym G=(G,\tp_{Sym})$. Obviously, $\tp_{Sym}$ is a topology that is stronger than $\tp$.

\begin{proposition}\label{p:apt:3}
Let $G$ be a group with a topology.
Then
\begin{enumerate}
\item
if $G$ is a semitopological group, then $\sym G$ is a quasitopological;
\item
if $G$ is a paratopological group, then $\sym G$ is topological;
\item
$\sym G$ is homeomorphic to $S_G$;
\item
if $G$ is an almost paratopological $T_1$ group, then
\begin{enumerate}
\item
$\sym G$ embeds closed in $G^2$;
\item
$\sym G$ is a Hausdorff quasitopological group.
\end{enumerate}
\end{enumerate}
\end{proposition}
\begin{proof}
Let $\tp$ be the topology of $G$ and $\tp_{Sym}$ be the topology of $\sym G$.

(1). If $g\in G$, $U,V\in\tp$ and $W=U\cap V^{-1}\in \tp_{Sym}$ then
\begin{align*}
gW&=gU\cap (Vg^{-1})^{-1}\in \tp_{Sym},
\\
Wg&=Ug\cap (g^{-1}V)^{-1}\in \tp_{Sym},
\\
W^{-1}&= V\cap U^{-1} \in \tp_{Sym}.
\end{align*}
Hence $\sym G$ is quasitopological.

(2). Follows from (1).

(3). Let's put
  \[
  \gj: \sym G \to S_G, \ x \mapsto (x,x^{-1}).
  \]
The mapping $\gj$ is a homeomorphism since
\[
\gj^{-1}(S_G\cap (U\times V)) = U\cap V^{-1}.
\]
for $U,V\in\tp$.

(4). From (3) and Proposition \ref{p:apt:2}(4) follows (a).
It follows from (1) that $G$ is a quasitopological group. Let us show that $\sym G$ is a Hausdorff space.
Let $e\neq g\in G$. Since $G$ is an almost paratopological $T_1$ group, then
$g\notin U^2$ for some $U\in\nne$. Let $S=U\cap U^{-1}$. Then $S$ is an open neighborhood of the identity in $\sym G$ and $eS\cap gS=\es$.
\end{proof}

\begin{proposition}\label{p:apt:4}
Let $G$ and $H$ be groups with topology and $\ph: G\to H$ be a topology-preserving homomorphism.
Let $\cP$ be one of the enumerated classes of groups: semitopological; quasitopological; paratopological; almost paratopological; compact. Then $G\in \cP$ if and only if $H\in\cP$.
\end{proposition}

\begin{theorem}\label{t:apt:2}
A compact almost paratopological group is a topological group.
\end{theorem}
\begin{proof}
Let $Q$ be a compact almost paratopological group. We set $H=\clp e$.
Theorem \ref{t:rsa:2} implies that $H$ is a normal closed antidiscrete subgroup,
the quotient group $G=Q/H$ is a $T_1$ compact semitopological group, and the quotient mapping $\ph: Q\to G$ is a topology-preserving homomorphism.
The Proposition \ref{p:apt:4} implies that $G$ is a compact almost paratopological $T_1$ group and $G$ is a topological group if and only if $Q$ is a topological group. Therefore, to prove the theorem, it suffices to prove that $G$ is a topological group.

It follows from the Proposition \ref{p:apt:3}(1) that $\sym G$ is quasitopological. The Proposition \ref{p:apt:3}(4)(a) implies that $\sym G$ embeds in $G^2$ in a closed manner, and hence the group $\sym G$ is compact. It follows from the Proposition \ref{p:apt:3}(4)(b) that $\sym G$ is Hausdorff. Hence $\sym G$ is a compact Hausdorff semitopological group. It follows from the Ellis theorem \cite[Theorem 2]{Ellis1957-2} that $\sym G$ is a topological group.

Let $\tp$ be the topology of $G$ and $\tp_{Sym}$ be the topology of $\sym G$.

Let us show that $G$ is a Hausdorff space.
Let $e\neq g\in G$.
Since $G$ is an almost paratopological $T_1$ group, it follows from Proposition \ref{p:apt:2}(4) that $E_G=\sset e$ and hence $g\notin \cl{U ^{-1}}$ for some $U\in\nne$.

\begin{claim}
$e\in\Int \cl{U^{-1}}$.
\end{claim}
\begin{proof}
Since $\sym G$ is a topological group and $U^{-1}\in\tp_{Sym}$, then $S^2\subset U^{-1}$ for some $S\in\tp_{ Sym}$ for which $e\in S=S^{-1}$. Since $\sym G$ is compact, then $G=\bigcup_{i=1}^n x_i S$ for some $x_1,x_2,...,x_n\in G$. A topological space cannot be the union of a finite number of nowhere dense sets, so $\Int \cl{x_i S}\neq\es$ for some $i$. Hence $\Int \cl{S}\neq\es$. Let $q\in S\cap \Int \cl{S}$. Then $e\in \Int\cl{q^{-1}S}\subset \cl{S^2}\subset \cl{U^{-1}}$. Hence $e\in \Int\cl{U^{-1}}$.
\end{proof}

Set $U_g=G\setminus \cl{U^{-1}}$ and $U_e=\Int \cl{U^{-1}}$. Then $g\in U_g\in\tp$, $e\in U_e\in\tp$ and $U_g\cap U_e=\es$.
Thus, the space $G$ is Hausdorff. The group $G$ is a Hausdorff compact semitopological group. It follows from Ellis \cite[Theorem 2]{Ellis1957-2} that $G$ is a topological group.
\end{proof}

Theorems \ref{t:apt:1} and \ref{t:apt:2} imply the following assertion.

\begin{cor}[{\cite[Lemma 6]{Ravsky2019}}]\label{c:apt:1}
A compact paratopological group is a topological group.
\end{cor}

\section{Spaces with separately continuous Mal'tsev operation}\label{sec:maltcev}
For the space $X$ we define the following property:
\begin{enumerate}
\item[\pUs/]
Let $\set{x_\al:\al<\om_1}\subset X$ and for $\al<\om_1$ let $\cF_\al$ be at most a countable family of closed subsets of $X$. Then there exists $\be<\om_1$ for which the following condition is satisfied:
\begin{enumerate}
\item[$\rm(*)$]
exists $y\in\cl{\set{x_\al:\al<\be}}$ such that
if $\ga<\be$, $F\in \cF_\ga$ and $x_\be\in F$, then $y\in F$.
\end{enumerate}
\end{enumerate}

The \pUs/ property is a strengthening of the $\rm (P)$ property from \cite{ReznichenkoUspenskij1998,uspenskij1989} in the class of Tikhonov spaces. The \pUs/ property can be used for regular spaces.

\begin{proposition}\label{p:maltcev:cel:1}
Let $X$ be a regular $\Si$-space with  countable extend. Then \pUs/ is true for $X$.
\end{proposition}
\begin{proof}
Let $\cS$ and $\cC$ be as in Statement \ref{a:maltcev:cel:1}.
We can assume that the family $\cS$ is closed under finite intersections.
Let $\set{x_\al:\al<\om_1}\subset X$ and for $\al<\om_1$ let $\cF_\al$ be at most a countable family of closed subsets of $X$. Denote
\begin{align*}
F^*_\ga &= \set{\bigcap\cF: \cF\subset \bigcup_{\al<\ga}\cF_\al,\ |\cF|<\om, \bigcap\cF\neq\es},
\\
X_\ga &= \set{x_\al: \al < \ga}
\end{align*}
for $\ga \leq \om_1$.
By induction on $n\in\om$ we construct an increasing sequence of countable ordinals $\sqnn \be$ such that for $n>0$ the following condition is satisfied:
\begin{enumerate}
\item[$(P_n)$]
if $x_\ga\in S\cap F$ for $\ga<\om_1$, $S\in\cS$ and $F\in \cF^*_{\be_{n-1}}$, then there exists $y\in\cl{X_{\be_{n}}}$ such that $y\in S\cap F$.
\end{enumerate}
Let $\be_0=\om$. Suppose that $n>0$ and $\be_0 < \be_1 < ... < \be_{n-1} < \om_1$ are constructed.
For $S\in\cS$ and $F\in \cF^*_{\be_{n-1}}$ we denote
\[
A_{S,F}=\set{\al<\om_1: x_\al\in S\cap F}.
\]
If $A_{S,F}\neq\es$, then we denote $\al_{S,F}=\min A_{S,F}$.
Let us put
\[
\be_n=\sup \set{\al_{S,F}: S\in\cS, F\in \cF^*_{\be_{n-1}}\text{ and }A_{S,F} \neq\es}+1.
\]

The sequence $\sqnn \be$ is constructed. Let $\be= \sup \set{\be_n: n\in\om}$.
Let us check $\rm(*)$ in \pUs/ definition.
There is $C\in\cC$ for which $x_\be\in C$.
Let us put
\begin{align*}
\cS'&=\set{S\in \cS: C\subset S}=\set{S_n':n\in\om},
\\
\cF'&=\set{F\in \cF^*_{\be}: x_\be\in F} =\set{F_n':n\in\om}.
\end{align*}

Let $n\in \om$.
Let us put
\begin{align*}
S_n&=\bigcap_{i=0}^n S_i',
&
F_n&=\bigcap_{i=0}^n F_i',
\\
\al_n&=\al_{S_n,F_n},
&
y_n &= x_{\al_n}.
\end{align*}
Since $\cF^*_{\be}=\bigcap_{m\in\om}\cF^*_{\be_m}$, then $F_n\in \cF^*_{\be_m}$ for some $m\in\om$.
Since $\be\in A_{S_n,F_n}\neq \es$, then $\al_n=\al_{S_n,F_n}\leq \be_{m+1}<\be$. Hence $y_n\in X_\be$.

It follows from the definition of the families $\cS$ and $\cC$ that the sequence $\sqnn y$ accumulates to some point $y\in C\cap\bigcap_{n\in\om} F_n$. Since $\sqnn y\subset X_\be$, then $y\in \cl{X_\be}$.

Let $F\in \cF_\ga$ for $\ga<\be$ and $x_\be\in F$.
Then $F\in \cF'$ and $F=F'_m\supset F_m$ for some $m\in\om$.
Hence,
\[
y\in \bigcap_{n\in\om} F_n \subset F_m\subset F.
\]
\end{proof}

\begin{proposition}\label{p:maltcev:cel:2}
Let $X$ be a regular space with a separately continuous Mal’tsev operation and let $X$ satisfy \pUs/.
Then $X$ is an $\om$-cellular space.
\end{proposition}
\begin{proof}
Let us assume the opposite. Then there is a family $\set{K_\al':\al<\om_1}$ of non-empty sets of type $G_\de$, such that
\[
K'_\be\not\subset \cl{\bigcup_{\al<\be}K'_\al}
\]
for $\be<\om_1$. Let us choose
\[
x_\be \in K'_\be\setminus \cl{\bigcup_{\al<\be}K'_\al}
\]
and a sequence $\sq{U_{\be,n}}{\nom}$ of open sets $X$, such that
\[
x_\be \in U_{n+1} \subset \cl{U_{n+1}}\subset U_n
\]
for $\nom$.
Then
\[
x_\be \in K_\be=\bigcap_\nom U_{\be,n}\subset K'_\be
\]
And
\begin{equation}\label{eq:maltcev:cel:0}
x_\be \notin \cl{\bigcup_{\al<\be}K_\al}.
\end{equation}
For $\al,\ga<\om_1$ we put
  \[
  h_{\al,\ga}: X\to X,\ x\mapsto M(x,x_\al,x_\ga).
  \]
Note that
\begin{align*}
h_{\al,\ga}(x_\al)&=x_\ga,
&
h_{\al,\al}=\id_X.
\end{align*}
Let us put
\begin{align*}
\cP_\be &= \set{ h_{\al,\ga}^{-1}(X\setminus{U_{\ga,n})}: \al,\ga<\be \text{ and } n<\om},
\\
\cF_\be &= \cP_{\be+1}.
\end{align*}
Since the condition \pUs/ is satisfied for $X$, then there exists $\be<\om_1$ and
\[
y\in\cl{\set{x_\al:\al<\be}},
\]
such that
if $\ga<\be$, $F\in \cF_\ga$ and $x_\be\in F$, then $y\in F$.
Then
\begin{equation}\label{eq:maltcev:cel:1}
\text{if } x_\be\in F\in \cP_\be,\text{ then }y\in F.
\end{equation}
Let us put
  \[
  y_\ga = M(x_\be,y,x_\ga)
  \]
for $\ga<\be$.
\begin{claim}
$y_\ga\in K_\ga$.
\end{claim}
\begin{proof}
Suppose $y_\ga\notin K_\ga$. Then $y_\ga\notin U_{\ga,n}$ for some $n\in\om$.
Let us put
\begin{align*}
U_1&=\set{x\in X: M(y,x,x_\ga)\in U_{\ga,n+1}},
\\
U_2&=\set{x\in X: M(x_\be,x,x_\ga)\in X\setminus \cl{U_{\ga,n+1}}}.
\end{align*}
The sets $U_1$ and $U_2$ are open. Because
\begin{align*} x_\ga=M(y,y,x_\ga) & \in U_{\ga,n+1}, 
\\
y_\ga=M(x_\be,y,x_\ga) & \notin U_{\ga,n}\supset \cl{U_{\ga,n+1}},
\end{align*}
then $y\in U_1\cap U_2$.
The set $U_1\cap U_2$ is an open neighborhood of $y$.
Since $y\in\cl{\set{x_\al:\al<\be}}$, then $x_\al\in U_1\cap U_2$ for some $\al<\be$.
Then
\begin{align*}
h_{\al,\ga}(y)&=M(y,x_\al,x_\ga)\in U_{\ga,n+1},
\\
h_{\al,\ga}(x_\be)&=M(x_\be,x_\al,x_\ga)\in X\setminus \cl{U_{\ga,n+1}}\subset X \setminus U_{\ga,n+1}.
\end{align*}
Hence
\begin{align*}
x_\be & \in h_{\al,\ga}^{-1}(X\setminus U_{\ga,n+1})=F\in \cP_\be,
\\
y & \in h_{\al,\ga}^{-1}({U_{\ga,n+1}})= X\setminus F.
\end{align*}
Contradiction with (\ref{eq:maltcev:cel:1}).
\end{proof}
Let us put
  \[
  h: X\to X, \ x \mapsto M(x_\be,y,x).
  \]
Since $h$ is continuous, $y_\ga=h(x_\ga)$ for $\ga<\be$ and $y\in \cl{\set{x_\ga: \ga<\be}}$, then
\[
x_\be=M(x_\be,y,y)=h(y)\in \cl{h(\set{x_\ga: \ga<\be})}=\cl{\set{y_\ga:\ga<\be}}.
\]
It follows from the claim that
\[
x_\be \in \cl{\bigcup_{\ga<\be}K_\ga}.
\]
Contradiction with (\ref{eq:maltcev:cel:0}).
\end{proof}
 Propositions \ref{p:maltcev:cel:1} and \ref{p:maltcev:cel:2} imply the following sentence.
\begin{theorem}\label{t:maltcev:cel:1}
Let $X$ be a regular $\Si$-space with countable extend and separately continuous Mal’tsev operation.
Then $X$ is an $\om$-cellular space.
\end{theorem}

\begin{cor}\label{c:maltcev:cel:1}
Let $X$ be a regular (closely $\si$-)countably compact space with separately continuous Mal’tsev operation.
Then $X$ is an $\om$-cellular space.
\end{cor}

\section{ccc in groups}\label{sec:ccc}
Since the Mal’tsev operation $M_G(x,y,z)=xy^{-1}z$ on a quasitopological group is separately continuous, then Theorem \ref{t:maltcev:cel:1} implies the following assertion.

\begin{theorem}\label{t:ccc:1}
Let $G$ be a regular $\Si$ quasitopological group with countable extend.
Then $G$ is an $\om$-cellular space.
\end{theorem}

\begin{cor}\label{c:ccc:1}
Let $G$ be a regular quasitopological group.
If any of the following conditions is satisfied for $G$, then $G$ is an $\om$-cellular space:
\begin{enumerate}
\item
$G$ is closely $\si$-countably compact;
\item
$G$ is a Lindelöf $\Si$-space;
\item
$G$ is $\si$-compact.
\end{enumerate}
\end{cor}

\begin{theorem}\label{t:ccc:2}
Let $G$ be a regular almost paratopological group and $G^2$ be a $\Si$-space with  countable extend.
Then $G$ is an $\om$-cellular space.
\end{theorem}
\begin{proof}
The Proposition \ref{p:apt:3} implies that $\sym G$ embeds closed in $G^2$ and is a quasitopological group. Hence $\sym G$ is a regular $\Si$ quasitopological group with countable extend. Theorem \ref{t:ccc:1} implies that $\sym G$ is $\om$-cellular. Since $G$ is a continuous image of $\sym G$, then $G$ is $\om$-cellular.
\end{proof}

Since the square of a Lindelöf $\Si$-space is a Lindelöf $\Si$-space, then
the Theorem \ref{t:ccc:2} implies the following assertion.

\begin{theorem}\label{t:ccc:3}
Let $G$ be a regular Lindelöf $\Si$ almost paratopological group.
Then $G$ is an $\om$-cellular space.
\end{theorem}

\begin{cor}\label{c:ccc:2}
Let $G$ be a regular almost paratopological group.
If any of the following conditions is satisfied for $G$, then $G$ is an $\om$-cellular space:
\begin{enumerate}
\item
$G^2$ is closely $\si$-countably compact;
\item
$G$ is a Lindelöf $\Si$-space;
\item
$G$ is $\si$-compact.
\end{enumerate}
\end{cor}

\begin{question}\label{q:ccc:1}
Let $G$ be a semitopological group.
Which of the following conditions imply that $G$ is an $\om$-cellular space ($G$ is ccc)?
\begin{enumerate}
\item
$G$ is a  $\si$-countably compact (regular) (almost paratopological, paratopological, quasitopological) group;
\item
$G$ is a closely $\si$-countably compact regular (almost paratopological, paratopological) group;
\item
$G$ is a (closely $\si$-)countably compact (almost paratopological, paratopological, quasitopological) group;
\item
$G$ is a $\si$-compact (almost paratopological, paratopological, quasitopological) group;
\item
$G$ is a Lindelöf $\Si$ group;
\item
$G$ is a $\si$-compact regular group.
\end{enumerate}
\end{question}

The author thanks the referee for useful comments.

\bibliographystyle{plainnat}


\bibliography{almost_paratopological}
\end{document}